\documentclass[12pt]{amsart}
\usepackage{amssymb,latexsym,verbatim}

\theoremstyle{plain}
\newtheorem{thm}{Theorem}[section]
\newtheorem{lem}[thm]{Lemma}
\newtheorem{cor}[thm]{Corollary}
\newtheorem{prop}[thm]{Proposition}

\theoremstyle{definition}
\newtheorem{defn}[thm]{Definition}

\theoremstyle{remark}

\newtheorem*{rem}{Remark}
\newtheorem*{conj}{Conjecture}
\newtheorem*{que}{Question}

\DeclareMathOperator{\Ann}{Ann}

\DeclareMathOperator{\Ext}{Ext}
\DeclareMathOperator{\Supp}{Supp}
\DeclareMathOperator{\V}{V}
\DeclareMathOperator{\Hom}{Hom}

\DeclareMathOperator{\Coker}{Coker}

\DeclareMathOperator{\Ass}{Ass}
\DeclareMathOperator{\Coass}{Coass}
\DeclareMathOperator{\Min}{Min}
\DeclareMathOperator{\Max}{Max}
\DeclareMathOperator{\lc}{H}
\DeclareMathOperator{\End}{End}

\DeclareMathOperator{\Spec}{Spec}
\DeclareMathOperator{\G}{\Gamma}
\DeclareMathOperator{\Att}{Att}
\DeclareMathOperator{\E}{E}

\newcommand{\lo}{\longrightarrow}
\newcommand{\fa}{\mathfrak{a}}

\newcommand{\fm}{\mathfrak{m}}
\newcommand{\fp}{\mathfrak{p}}
\newcommand{\fq}{\mathfrak{q}}

\begin{document}

\title[local cohomology: cofiniteness, coassociated primes ]
{Cofiniteness and coassociated primes of local cohomology modules \\}

\author{Moharram Aghapournahr}
\address{ Moharram Aghapournahr\\Arak University\\Beheshti St, P.O. Box:879,
 Arak, Iran}
\email{m-aghapour@araku.ac.ir} \email{m.aghapour@gmail.com}

\author{Leif Melkersson}
\address{Leif Melkersson\\Department of Mathematics\\
         Link\"{o}ping University\\
         SE--581 83 Link\"{o}ping, Sweden}
\email{lemel@mai.liu.se}




\begin{abstract}
 Let $R$ be a noetherian ring,
$\fa$ an ideal of $R$ such that $\dim R/\fa=1$ and $M$ a finite $R$--module.
We will study cofiniteness and some other properties of
the local cohomology modules $\lc^{i}_{\fa}(M)$.
For an arbitrary ideal $\fa$ and an $R$--module $M$ (not necessarily finite),
we will characterize $\fa$--cofinite artinian local cohomology modules.
Certain sets of coassociated primes of top local cohomology modules over
local rings are characterized.
\end{abstract}

\maketitle

\section{Introduction}
Throughout $R$ is a commutative noetherian ring.
By a finite module we mean a finitely generated module.
For basic facts about commutative algebra see \cite{BH} and \cite{Mat}
and for local cohomology we refer to \cite{BSh}.

Grothendieck \cite{SGA2}, made the following conjecture:

\begin{conj}
For every ideal $\fa$ and every finite $R$--module $M$,
the module $\Hom_{R}(R/\fa,\lc^{n}_{\fa}(M))$ is finite for all $n$.
\end{conj}

Hartshorne \cite{Ha} showed that this is false in general.
However, he defined an $R$--module $M$ to be {\it $\fa$--cofinite} if
$\Supp_R(M)\subset \V{(\fa)}$ and $\Ext^i_{R}(R/\fa,M)$
is finite (finitely generated) for each $i$
and he asked the following question:

\begin{que} If $\fa$ is an ideal of $R$ and
$M$ is a finite $R$--module.
When is $\Ext^i_{R}(R/\fa,\lc^{j}_{\fa}(M))$ finite for every $i$ and $j$ ?
\end{que}

Hartshorne \cite{Ha} showed that if
$(R,\fm)$ is a complete regular local ring and $M$ a finite $R$--module, then
$\lc^{i}_{\fa}(M)$ is $\fa$--cofinite in two cases:

\begin{flushleft}
(a) If $\fa$ is a nonzero principal ideal, and
\end{flushleft}
(b) If $\fa$ is a prime ideal with $\dim R/\fa=1$.

Yoshida \cite{Yo} and Delfino and Marley \cite{DM} extended (b) to all
dimension one ideals $\fa$ of an arbitrary local ring $R$.

In \ref{C:cofmin2}, we give a characterization of the
$\fa$--cofiniteness of these local cohomology modules
when $\fa$ is a one-dimensional ideal in a non-local ring.
In this situation we also prove
in  \ref{T:Hinz},  that these local cohohomology modules
always belong to a class introduced by Z\"{o}schinger in \cite{Zrmm}.

Our main result in this paper is \ref{T:artcof}, where we
for an arbitrary ideal $\fa$ and an $R$--module $M$ (not necessarily finite),
 characterize the artinian $\fa$--cofinite local cohomology modules
(in the range $i< n$).
With the additional assumption that $M$ is finitely generated,
the characterization is also given by the existence of certain filter-
regular sequences.

The second author has in \cite[Theorem 5.5]{Mel} previously
characterized artinian local cohomology modules, (in the same
range). In case the module $M$ is not supposed to be finite, the
two notions differ. For example let $\fa$ be an ideal of a local
ring $R$, such that $\dim(R/\fa)>0$ and let
  $M$ be the injective hull of the residue field of $R$.
The module $\lc^0_\fa(M)$, which is equal to $M$, is artinian.
However it is not $\fa$--cofinite, since
$0\underset{M}:\fa$ does not have finite length.

An $R$--module $M$ has {\it finite Goldie dimension} if
$M$ contains no infinite direct sum of submodules.
For a commutative noetherian ring this can be expressed in two other ways,
namely that the injective hull
$\E(M)$ of $M$ decomposes as a finite direct sum of indecomposable
injective modules or that $M$ is an essential extension of a finite submodule.

A prime ideal $\fp$ is said to be {\it coassociated} to $M$ if
$\fp=\Ann_R({M/N})$ for some $N\subset M$ such that
$M/N$ is artinian and is said to be {\it attached} to $M$ if
$\fp=\Ann_R({M/N})$ for some arbitrary submodule $N$ of $M$, equivalently
$\fp=\Ann_R({M/{\fp}M})$. The set of these prime ideals are denoted by
$\Coass_R(M)$ and $\Att_R(M)$ respectively. Thus
$\Coass_R(M)\subset \Att_R(M)$ and the two sets are equal when
$M$ is an artinian module.
The two sets behave well with respect to exact sequences.
If $0\rightarrow M^{\prime}\rightarrow M\rightarrow M^{\prime\prime}
\rightarrow 0$ is an exact sequence, then
$$\Coass_R(M^{\prime\prime})\subset \Coass_R(M)\subset{\Coass_R(M^{\prime})\cup \Coass_R(M^{\prime\prime})}$$
\begin{center}
and
\end{center}
\begin{center}
$\Att_R(M^{\prime\prime})\subset \Att_R(M)\subset{\Att_R(M^{\prime})
\cup \Att_R(M^{\prime\prime})}.$
\end{center}
There are equalities $\Coass_R(M\otimes_R{N})= \Coass_R(M)\cap \Supp_R(N)$ and
$\Att_R(M\otimes_R{N})= \Att_R(M)\cap \Supp_R(N)$,
whenever the module $N$ is required to be finite.
We prove the second equality in \ref{L:att}.
In particular $\Coass_R(M/{\fa}M)=\Coass_R(M)\cap \V(\fa)$ and
$\Att_R(M/{\fa}M)=\Att_R(M)\cap \V(\fa)$ for every ideal $\fa$.
Coassociated and attached prime ideals have been studied in particular by
Z\"{o}schinger, \cite{Zrkoass} and \cite{Zrlk}.

In \ref{C:attcoass} we give a characterization of certain sets of
coassociated primes of the highest nonvanishing local cohomology
module $\lc_{\fa}^t(M)$, where $M$ is a finitely generated module
over a complete local ring. In case it happens that $t=\dim M$,
the characterization is given in \cite[Lemma 3]{DM}. In that case
the top local cohomology module is always artinian, but in
general the top local cohomology module  is not artinian
 if   $t<\dim M$.


\section{Main results}

First we extend a result by Z\"{o}schinger \cite[Lemma 1.3]{Zrko}
with a much weaker condition. Our method of proof is also quite different.

\begin{prop}\label{P:fgwl}

Let $M$ be a module over the noetherian ring $R$.
The following statements are equivalent:
\begin{enumerate}
  \item[(i)]$M$ is a finite $R$--module.
  \item[(ii)]$M_\fm$ is a finite $R_\fm$--module for all $\fm{\in}\Max{R}$
and\\$\Min_R(M/N)$ is a finite set for all finite submodules $N\subset M$.
\end{enumerate}
\end{prop}
\begin{proof}
The only nontrivial part is (ii)$\Rightarrow$ (i).

Let $\mathcal F$ be the set of finite submodules of $M$.
For each $N \in\mathcal F$ the set $\Supp_R(M/N)$ is closed in $\Spec (R)$,
since $\Min_R(M/N)$ is a finite set.
Also it follows from the hypothesis that, for each
$\fp\in \Spec (R)$ there is $N\in \mathcal F$ such that $M_{\fp}=N_{\fp}$,
that is $\fp\notin \Supp_R(M/N)$.
This means that ${\bigcap}_{N\in \mathcal F}{\Supp_R(M/N)}=\varnothing$.
Now $\Spec (R)$ is a quasi-compact topological space.
Consequently $\bigcap_{i=1}^r\Supp_R(M/N_i)=\varnothing$ for some
$N_1,...,N_r\in \mathcal F$.
We claim that $M=N$, where $N=\sum_{i=1}^r N_i$.
Just observe that $\Supp_R(M/N)\subset \Supp_R(M/N_i)$ for each $i$,
and therefore $\Supp_R(M/N)=\varnothing$.
\end{proof}

\begin{cor}\label{C:cofmin1}
Let $M$ be an $R$--module such that
$\Supp M\subset \V(\fa)$ and $M_\fm$ is
${\fa}R_{\fm}$--cofinite for each maximal ideal $\fm$.
The following statements are equivalent:
\begin{enumerate}
  \item[(i)] $M$ is $\fa$--cofinite.
  \item[(ii)] For all $j$,  $\Min_R(\Ext^{j}_{R}(R/\fa,M)/T)$ is a
finite set for each finite submodule $T$ of $\Ext^{j}_{R}(R/\fa,M)$.
\end{enumerate}
\end{cor}
\begin{proof}
The only nontrivial part is (ii)$\Rightarrow$ (i).

Suppose $\fm$ is a maximal ideal of $R$.
By hypothesis $M_\fm$ is ${\fa}R_{\fm}$--cofinite.
Therefore $\Ext^{j}_{R}(R/\fa, M)_{\fm}$ is a finite
$R_\fm$--module for all $j$.
Hence by \ref{P:fgwl} \  $\Ext^{j}_{R}(R/\fa,M)$ is finite for all $j$.
Thus $M$ is $\fa$--cofinite.
\end{proof}

\begin{cor}\label{C:cofmin2}
Let $\fa$ an ideal of $R$ such that $\dim R/\fa=1$,
$M$ a finite $R$--module and $i\geq 0$.
The following statements are equivalent:
\begin{enumerate}
  \item[(i)] $\lc^{i}_{\fa}(M)$ is $\fa$--cofinite.
  \item[(ii)] For all $j$,
$\Min_R(\Ext^{j}_{R}(R/\fa, \lc^{i}_{\fa}(M))/T)$
is a finite set for each finite submodule $T$ of
$\Ext^{j}_{R}(R/\fa, \lc^{i}_{\fa}(M))$.
\end{enumerate}
\end{cor}

\begin{proof}
For all maximal ideals
$\fm$, $\lc^{i}_{\fa}(M)_\fm \cong \lc^{i}_{{\fa}R_{\fm}}(M_\fm)$.
By \cite[Theorem 1]{DM} $\lc^{i}_{{\fa}R_{\fm}}(M_\fm)$ is
${\fa}R_{\fm}$--cofinite.
\end{proof}

A module $M$ is {\it weakly Laskerian}, when for each submodule $N$ of $M$
the quotient $M/N$ has just finitely many associated primes, see \cite{DiM}.
A module $M$ is {\it $\fa$--weakly cofinite} if
$\Supp_R(M)\subset \V(\fa)$ and
$\Ext^{i}_{R}(R/\fa, M)$ is weakly Laskerian for all $i$.
Clearly each $\fa$--cofinite module is
$\fa$--weakly cofinite but the converse is not true in general
see \cite[Example 3.5 (i) and (ii)]{DiM2}.

\begin{cor}\label{c:wlcof}
If $\lc^{i}_{\fa}(M)$ {\rm (}with $\dim R/\fa=1${\rm )} is an
$\fa$--weakly cofinite module, then it is also $\fa$--cofinite.
\end{cor}
Next we will introduce a subcategory of the category of $R$--modules
that has been studied by Z\"{o}schinger in \cite[Satz 1.6]{Zrmm}.

\begin{thm}\label{T:classz}\textbf{{\rm (}Z\"{o}schinger{\rm )}}
For any $R$--module $M$ the following are equivalent:
\begin{enumerate}
  \item[(i)] $M$ satisfies the minimal condition for submodules $N$ such that
$M/N$ is soclefree.
  \item[(ii)]For any descending chain
$N_1\supset N_2\supset N_3\supset \dots$ of submodules of $M$,
there is $n$ such that the quotients
$N_{i}/N_{i+1}$ have support in $\Max R$ for all $i\geq n$.
  \item[(iii)] With $L(M)=\underset{\fm\in \Max R}{\bigoplus}\G_{\fm}(M)$,
the module $M/L(M)$ has finite Goldie dimension,
and $\dim R/\fp \leq 1$ for all $ \fp \in \Ass_R(M)$.
\end{enumerate}

If they are fulfilled, then for each monomorphism $f:M\lo M$,
$$ \Supp_R(\Coker f)\subset \Max R.$$
\end{thm}

We will say that $M$ is in the class
$\mathcal Z$ if $M$ satisfies the equivalent conditions in \ref{T:classz}.

A module {M} is {\it soclefree} if it has no simple submodules,
or in other terms $\Ass M\cap\Max R=\varnothing$.
For example if $M$ is a module over the local ring $(R,\fm)$
then the module $M/{\G_{\fm}(M)}$, where $\G_{\fm}(M)$ is the submodule of
$M$ consisting of all elements of $M$ annihilated by some high power
${\fm}^n$ of the maximal ideal $\fm$, is always soclefree.

\begin{prop}\label{T:serrez}
The class $\mathcal Z$ is a Serre subcategory of the category of $R$--modules,
that is $\mathcal Z$ is closed under
taking submodules, quotients and extensions.
\end{prop}
\begin{proof}
The only difficult part is to show that $\mathcal Z$ is closed
under taking extensions. To this end let $0\lo M^{\prime}\overset
f\lo M\overset g\lo M^{\prime\prime}\lo 0$ be an exact sequence
with $M^{\prime},M^{\prime\prime}\in\mathcal Z$ and let
$N_1\supset N_2\supset ...$ be a descending chain of submodules
of $M$. Consider the descending chains $f^{-1}(N_1)\supset
f^{-1}(N_2)\supset ...$ and $g(N_1)\supset g(N_2)\supset ...$ of
submodules of $M^{\prime}$ and $M^{\prime\prime}$ respectively. By
(ii) there is $n$ such that
$\Supp_R(f^{-1}(N_i)/f^{-1}(N_{i+1}))\subset \Max R$ and
$\Supp_R(g(N_i)/g(N_{i+1}))\subset \Max R$ for all $i\geq n$. We
use the exact sequence
$$0\lo f^{-1}(N_i)/f^{-1}(N_{i+1})\lo N_i/N_{i+1}\lo g(N_i)/g(N_{i+1})\lo 0.$$
to conclude that $\Supp_R(N_i/N_{i+1})\subset \Max R$ for all $i\geq n$.
\end{proof}

\begin{thm}\label{T:Hinz}
Let $N$ be a module over a noetherian ring
$R$ and $\fa$ an ideal of $R$ such that $\dim{R/\fa}=1$.
If $N_\fm$ is ${\fa}R_\fm$--cofinite for all $\fm\in \Max R$,
then $N$ is in the class $\mathcal Z$.
In particular, if $M$ is a finite $R$--module then
$\lc^{i}_{\fa}(M)$ is in the class $\mathcal Z$ for all $i$.
\end{thm}

\begin{proof}
Let $X=N/L(N)$. Note that
$\Ass_R(X)\subset \Min{\fa}$ and therefore is a finite set. Since
$$
\E(X)=\underset{\fp\in \Ass_R(X)}\bigoplus \E(R/\fp)^{\mu^{i}(\fp,X)},
$$
 it is enough to prove that $\mu^{i}(\fp,X)$ is finite for all
$\fp\in \Ass_R(X)$.
This is clear, since each $\fp\in \Ass_R(X)$ is minimal over $\fa$ and
therefore $X_\fp \cong N_{\fp}$ which is,
${\fa}R_{\fp}$--cofinite, i.e. artinian over $R_{\fp}$.
\end{proof}

Given  elements $x_1,\dots,x_r$ in $R$, we denote by
$\lc^{i}(x_1,\dots,x_r;M)$  the $i$'th  Koszul cohomology module  of
the $R$--module $M$.
The following lemma is used in the proof of \ref{T:artcof}.

\begin{lem}\label{L:inj}

Let $E$ be an injective module.
If $\lc^0(x_1,\dots,x_r ; E)=0$, then $\lc^i(x_1,\dots,x_r ; E)=0$ for all $i$.
\end{lem}
\begin{proof}
We may assume that $E=\E(R/\fp)$ for some prime ideal $\fp$,
since $E$ is a direct sum of modules of this form, and
Koszul cohomology preserves (arbitrary) direct sums.

Put $\fa=(x_1,\dots,x_r)$.
By hypothesis $0:_E{\fa}=0$, which means that
$\fa\not\subset \fp$. Take an element $s\in \fa\setminus \fp$.
It acts bijectively on $E$,
hence also on $\lc^i(x_1,\dots,x_r ; E)$ for each $i$.
But $\fa\subset \Ann_R({\lc^i(x_1,\dots,x_r ; E)})$ for all $i$,
so the element $s$ therefore acts as the zero homomorphism on each
$\lc^i(x_1,\dots,x_r ; E)$. The conclusion follows.
\end{proof}

First we state the definition, given in \cite{Mel},
of the notion of filter regularity on modules (not necessarily finite)
over any noetherian ring.
When $(R,\fm)$ is local and $M$ is finite,
it yields the ordinary notion of filter-regularity, see \cite{CST}.

\begin{defn}
Let $M$ be a module over the noetherian ring $R$.
An element $x$ of $R$ is called filter-regular on $M$
if the module $0:_M{x}$ has finite length.

A sequence $x_1,...,x_s$ is said to be filter regular on $M$ if
$x_j$ is filter-regular on $M/(x_1,...,x_{j-1})M$ for $j=1,...,s$.
\end{defn}

The following theorem yields a characterization of
artinian cofinite local cohomology modules.

\begin{thm}\label{T:artcof}
Let $\fa=(x_1,...,x_r)$ be an ideal of a noetherian ring $R$ and let
$n$ be a positive integer.
For each $R$--module $M$ the following conditions are equivalent:
\begin{enumerate}
  \item[(i)] $\lc^{i}_{\fa}(M)$ is artinian and $\fa$--cofinite for all $i<n$.
  \item[(ii)]$\Ext^{i}_{R}(R/\fa,M)$ has finite length for all  $i<n$.
  \item[(iii)]The Koszul cohomology modules $\lc^i(x_1,\dots,x_r ; M)$
    has finite length for all  $i<n$.
\end{enumerate}
When $M$ is finite these conditions are also equivalent to:
\begin{enumerate}
  \item[(iv)] $\lc^{i}_{\fa}(M)$ is artinian for all $i<n$.
  \item[(v)] There is a sequence of length $n$ in $\fa$ that is
             filter-regular on $M$.
\end{enumerate}
\end{thm}

\begin{proof}
We use induction on $n$.
When $n=1$ the conditions (ii) and (iii) both say that
$0:_M{\fa}$ has finite length,
and they are therefore equivalent to (i) \cite[Proposition 4.1]{Mel}.

Let $n> 1$ and assume that the conditions are equivalent when
$n$ is replaced by $n-1$.
Put $L=\G_{\fa}(M)$ and $\overline{M}=M/L$ and form the exact sequence
$0\lo L\lo M\lo \overline{M}\lo 0$.
We have $\G_{\fa}(\overline{M})=0$ and
$\lc^{i}_{\fa}(\overline{M})\cong \lc^{i}_{\fa}(M)$ for all $i> 0$.
There are exact sequences
$$
\Ext^{i}_{R}(R/\fa,L)\rightarrow \Ext^{i}_{R}(R/\fa,M)\rightarrow
\Ext^{i}_{R}(R/\fa,\overline{M})\rightarrow \Ext^{i+1}_{R}(R/\fa,L)
$$
\begin{center}
and
\end{center}
\begin{center}
$\lc^i(x_1,\dots,x_r ; L)\rightarrow \lc^i(x_1,\dots,x_r ; M)
\rightarrow \lc^i(x_1,\dots,x_r ; \overline{M}) \rightarrow
\lc^{i+1}(x_1,\dots,x_r ; L)$
\end{center}
Because $L$ is artinian and $\fa$--cofinite the outer terms of both
exact sequences have finite length.
Hence $M$ satisfies one of the conditions if and only if
$\overline{M}$ satisfies the same condition.
We may therefore assume that $\G_{\fa}(M)=0$.

Let $E$ be the injective hull of $M$ and put $N=E/M$.
Consider the exact sequence
$0\lo M\lo E\lo N\lo 0$.
We know that $0:_M{\fa}=0$.
Therefore $0:_E{\fa}=0$ and $\G_{\fa}(E)=0$.
Consequently there are isomorphisms for all $i\geq 0$:
$$
\lc^{i+1}_{\fa}(M)\cong \lc^{i}_{\fa}(N),
$$

$$
\Ext^{i+1}_{R}(R/\fa,M)\cong \Ext^{i}_{R}(R/\fa,N)
$$
\begin{center}
and
\end{center}
\begin{center}
$\lc^{i+1}(x_1,\dots,x_r ; M)\cong \lc^i(x_1,\dots,x_r ; N).$
\end{center}
In order to get the third isomorphism, we used that
$\lc^i(x_1,\dots,x_r ; E)=0$ for all $i\geq 0$ (\ref{L:inj}). Hence
$M$ satisfies one of the three conditions if and only if $N$
satisfies the same condition, with $n$ replaced by $n-1$. By
induction, we may therefore conclude that the module $M$ satisfies
all three conditions if it satisfies one of them.

Let now $M$ be a finite module.

(ii)$\Leftrightarrow $(iv) Use \cite[Theorem 5.5 (i)
$\Leftrightarrow $(ii)]{Mel}.

(v)$\Rightarrow $(i) Use \cite[Theorem 6.4]{Mel}.

(i)$\Rightarrow $(v) We give a proof by induction on $n$. Put
$L=\G_{\fa}(M)$ and $\overline{M}=M/L$. Then $\Ass_R L=\Ass_R
M\cap\V(\fa)$ and $\Ass_R\overline M=\Ass_R M\setminus\V(\fa)$.
The module $L$ has finite length and therefore $\Ass_R L\subset\Max R$.
By prime avoidance take an element
$y_1\in\fa\setminus\bigcup_{\fp\in\Ass_R(\overline M)}{\fp}$. Then
$\Ass_R(0:_M{y_1})=\Ass_R(M)\cap\V(y_1)
= (\Ass_R L\cap\V(y_1))\cup(\Ass_R\overline M \cap\V(y_1)) \subset\Max R$,
Hence $0:_M{y_1}$ has finite length,
so the element $y_1\in\fa$
is filter regular on $M$.

Suppose $n> 1$ and take $y_1$ as above.

Note that $\lc^{i}_{\fa}(M)\cong \lc^{i}_{\fa}(\overline{M})$ for all
$i\geq 1$.
Thus we may replace $M$ by $\overline{M}$,
\cite[Proposition 6.3 (b)]{Mel},
and we may assume that $y_1$ is a non-zerodivisor on $M$.

The exact sequence
$0\rightarrow M\overset{y_1}\rightarrow M\rightarrow M/{y_1}M\rightarrow 0$
yields the long exact sequence

$$
\dots\lo \lc^{i-1}_{\fa}(M)\lo \lc^{i-1}_{\fa}(M/{y_1}M)
\lo \lc^{i}_{\fa}(M)\lo \dots.
$$

Hence $\lc^{i}_{\fa}(M/{y_1}M)$ is $\fa$--cofinite and artinian for
all $i< n-1$, by \cite[Corollary 1.7]{LMcof}. Therefore by the
induction hypothesis there exists $y_2,\dots,y_n$ in $\fa$, which is
filter-regular on $M/{y_1}M$. Thus $y_1,\dots,y_n$ is filter-regular
on $M$.
\end{proof}

\begin{rem}\label{r:cofart2} In \cite{AMel} we studied
the kernel and cokernel of the natural homomorphism
$f:\Ext_R^n(R/\fa ,M)\to\Hom_R(R/\fa ,\lc^n_{\fa}(M))$.
Applying the criterion  of \ref{T:artcof} we get that if
 $\Ext^{t-j}_{R}(R/\fa, \lc^{j}_{\fa}(M))$ has finite length for
$t=n,n+1$ and for all $j<n$,
then $\Ext^n_{R}(R/\fa,M)$ has finite length if and only if
$\lc^{n}_{\fa}(M)$ is $\fa$--cofinite artinian.
\end{rem}

Next we will study attached and coassociated prime ideals for the last
nonvanishing local cohomology module.
First we prove a lemma used in  \ref{C:attcoass}

\begin{lem}\label{L:att}
For all $R$--modules $M$  and for every finite $R$--module $N$,

$$
\Att_R(M\otimes_R{N})=\Att_R(M)\cap \Supp_R(N).
$$

\end{lem}

\begin{proof}
Let $\fp\in \Att_R(M\otimes_R{N})$, so
$ \fp=\Ann_R((M\otimes_R{N})\otimes_R{R/\fp})$.
However this ideal contains  both
$\Ann_R (M/{\fp M})$ and $\Ann_R (N)$
 and therefore
$\fp= \Ann_R(M/{\fp M})$ and $\fp\in\Supp_R(N)$.

Conversely let $\fp\in \Att_R(M)\cap \Supp_R(N)$.
Then $\fp=\Ann M/{\fp}M$\ and we want to show that
$ \fp=\Ann_R((M\otimes_R{N})\otimes_R{R/\fp})$. Since
\begin{center}
$(M\otimes_R{N})\otimes_R{R/\fp}\cong M/{\fp}M\otimes_{R/\fp}{N/{\fp}N}$,
\end{center}
we may assume that $R$ is a domain and $\fp=(0)$. Let $K$ be the
field of fractions of $R$. Then $\Ann M=0$ and $N\otimes_R{K}\neq
0$. Therefore the natural homomorphism $f: R\lo \End_R(M)$ is
injective and we have the following exact sequence
$$
0\lo \Hom_R(N,R)\lo \Hom_R(N,\End_R(M)).
$$
But $\Hom_R(N,\End_R(M))\cong \Hom_R(M\otimes_R{N},M)$.
Hence we get
\begin{center}

$\Ann_R({M\otimes_R{N}})\subset \Ann_R{\Hom_R(M\otimes_R{N},M)}
\subset \Ann_R{\Hom_R(N,R)}\subset
\Ann_R({\Hom_R(N,R)\otimes_R{K}}).$

\end{center}
On the other hand $\Hom_R(N,R)\otimes_R{K}\cong
\Hom_R(N\otimes_R{K},K)$, which is a nonzero vector space over $K$.
Consequently $\Ann_R({M\otimes_R{N}})=0$.
\end{proof}

\begin{thm}\label{T:attcoass}
Let $(R,\fm)$ be a complete local ring and let $\fa$ be an ideal of $R$.
Let $t$ be a nonnegative integer such
that $\lc^{i}_{\fa}(R)=0$ for all $i>t$.
\begin{enumerate}
  \item[(a)] If $\fp\in \Att_R(\lc^{t}_{\fa}(R))$ then $\dim R/{\fp}\geq t.$
  \item[(b)] If $\fp$ is a prime ideal such that $\dim R/{\fp}=t$,
    then the following conditions are equivalent:
\begin{enumerate}
  \item[(i)] $\fp\in \Coass_R(\lc^{t}_{\fa}(R))$.
  \item[(ii)] $\fp\in \Att_R(\lc^{t}_{\fa}(R))$.
  \item[(iii)] $\lc^{t}_{\fa}(R/\fp)\neq 0$.
  \item[(iv)] $\sqrt{\fa+\fp}=\fm$.
\end{enumerate}
\end{enumerate}
\end{thm}

\begin{proof}
(a) By the right exactness of the functor $\lc^{t}_{\fa}(-)$ we have
\begin{equation}\label{E:iso}
\lc^{t}_{\fa}(R/\fp)\cong \lc^{t}_{\fa}(R)/{\fp}\lc^{t}_{\fa}(R)
\end{equation}
If $\fp\in \Att_R(\lc^{t}_{\fa}(R))$, then
$\lc^{t}_{\fa}(R)/{\fp}\lc^{t}_{\fa}(R)\neq 0$.
Hence $\lc^{t}_{\fa}(R/\fp)\neq 0$ and
$\dim R/{\fp}\geq t.$

(b) Since $R/\fp$ is a complete local domain of dimension $t$,
the equivalence of (iii) and (iv) follows from the local
Lichtenbaum Hartshorne vanishing theorem.

If $\lc^{t}_{\fa}(R/\fp)\neq 0$, then by (\ref{E:iso})
$\lc^{t}_{\fa}(R)/{\fp}\lc^{t}_{\fa}(R)\neq 0$.
Therefore $\fp\subset \fq$ for some
$\fq\in \Coass_R(\lc^{t}_{\fa}(R))\subset \Att_R(\lc^{t}_{\fa}(R))$.
By (a) $\dim R/{\fq}\geq t= \dim R/{\fp}$, so we must have $\fp=\fq$.
Thus (iii) implies (i) and since always
$\Coass_R(\lc^{t}_{\fa}(R))\subset \Att_R(\lc^{t}_{\fa}(R))$, (i) implies (ii).

If (ii) holds then the module
$\lc^{t}_{\fa}(R)/{\fp}\lc^{t}_{\fa}(R)\neq 0$, since its annihilator is zero.
 Hence, using again the isomorphism (\ref{E:iso}), (ii) implies (iii).
\end{proof}

\begin{cor}\label{C:attcoass}
Let $(R,\fm)$ be a complete local ring, $\fa$ an ideal of $R$ and
$M$ a finite $R$--module and $t$ a nonnegative integer such that
$\lc^{i}_{\fa}(M)=0$ for all $i>t$.
\begin{enumerate}
  \item[(a)] If $\fp\in \Att_R(\lc^{t}_{\fa}(M))$ then $\dim R/{\fp}\geq t.$
  \item[(b)] If $\fp$ is a prime ideal in
    $\Supp_R(M)$ such that $\dim R/{\fp}=t$, then the following
conditions are equivalent:
\begin{enumerate}
  \item[(i)] $\fp\in \Coass_R(\lc^{t}_{\fa}(M))$.
  \item[(ii)] $\fp\in \Att_R(\lc^{t}_{\fa}(M))$.
  \item[(iii)] $\lc^{t}_{\fa}(R/\fp)\neq 0$.
  \item[(iv)] $\sqrt{\fa+\fp}=\fm$.
\end{enumerate}
\end{enumerate}
\end{cor}
\begin{proof}
Passing from $R$ to $R/\Ann M$, we may assume that $\Ann M=0$ and
therefore using Gruson's theorem, see \cite[Theorem 4.1]{V},
$\lc^{i}_{\fa}(N)=0$ for all $i>t$ and every $R$--module $N$. Hence
the functor $\lc^{t}_{\fa}(-)$ is right exact and therefore, since
it preserves direct limits, we get
$$
\lc^{t}_{\fa}(M)\cong M\otimes_R{\lc^{t}_{\fa}(R)}.
$$
The claims follow from \ref{T:attcoass} using the following equalities
$$
\Coass_R(\lc^{t}_{\fa}(M))=\Coass_R(\lc^{t}_{\fa}(R))\cap \Supp_R(M)
$$
 by \cite[Folgerung 3.2]{Zrmm} and
$$
\Att_R(\lc^{t}_{\fa}(M))=\Att_R(\lc^{t}_{\fa}(R))\cap \Supp_R(M)
$$
 by \ref{L:att}.
\end{proof}


\end{document}